\documentclass[12pt,a4paper]{article}
\usepackage{amssymb}

\usepackage{graphicx,amssymb,amsfonts,epsfig,amsthm,a4,amsmath,url}
\usepackage[latin1]{inputenc}

\newtheorem{thm}{Theorem}[section]
\newtheorem{cor}[thm]{Corollary}
\newtheorem{lem}[thm]{Lemma}
\newtheorem{clai}[thm]{Claim}
\newtheorem{prop}[thm]{Proposition}
\theoremstyle{definition}
\newtheorem{defn}[thm]{Definition}

\newtheorem{nota}[thm]{Notation}

\theoremstyle{remark}
\newtheorem{rem}[thm]{Remark}

\numberwithin{equation}{section}

\newcommand{\Z}{\mathbf{Z}}

\newcommand{\A}{\mathcal{A}}
\newcommand{\CC}{\mathcal{C}}
\newcommand{\N}{\mathbf{N}}
\newcommand{\R}{\mathbf{R}}
\newcommand{\C}{\mathbf{C}}

\newcommand{\bpr}{\noindent \textbf{Proof}: ~}

\newcommand{\epr}{~$\blacksquare$}

\newcommand{\eps}{\varepsilon}

\newcommand{\supp}{\textnormal{Supp}}

\title{Left inverses of matrices with polynomial decay.}
\author{Romain Tessera\footnote{This work was conducted in June 2007, while the author was visiting
the Bernoulli center in Lausanne. The author is supported by the
NSF grant DMS-0706486.} }
\date{\today}

\begin{document}

\baselineskip=16pt

\maketitle

\begin{abstract}
It is known that the algebra of Schur operators on $\ell^2$ (namely operators bounded on both $\ell^1$ and $\ell^{\infty}$) is not inverse-closed.
When $\ell^2=\ell^2(X)$ where $X$ is a metric space, one can consider elements of the Schur algebra with certain decay at infinity. For instance if $X$ has the doubling property, then Q. Sun has proved that
 the weighted Schur algebra $\A_{\omega}(X)$ for a strictly polynomial weight $\omega$ is  inverse-closed.
In this paper, we prove a sharp result on left-invertibility of the these operators. Namely, if an operator $A\in \A_{\omega}(X)$ satisfies
$$\|Af\|_p\succeq \|f\|_p$$ for some $1\leq p\leq \infty$, then it admits a left-inverse in $\A_{\omega}(X).$  The main difficulty here is to obtain the above inequality in $\ell^2$. The author was both motivated and
inspired by a previous work of Aldroubi, Baskarov and Krishtal \cite{ABK}, where similar results were obtained through different methods for $X=\Z^d$, under additional conditions on the decay. 
\end{abstract}




\section{Introduction}
In this paper, we study the left-invertibility of certain classes
of bounded linear operators $A:\ell^p(X)\to \ell^p(Y)$ where $X$
is a metric space and $Y$ is any set. 

We say that such an operator is bounded below in $\ell^p$ if
$$\lambda_p(A):=\inf_{f\neq 0}\frac{\|Af\|_p}{\|f\|_p}>0.$$
If $A$ is left invertible in $\ell^p$, i.e. if there exists a bounded linear map $B:\ell^p(Y)\to \ell^p(X)$ such that $BA=I$, then $A$ is clearly bounded below in $\ell^p$. But unless $p=2$, the converse is not true in general. Our main concern in this article will be to prove the converse in certain situations, namely when the matrix satisfies some {\it decay} condition. The first results of this kind were obtained in \cite{ABK}.
This type of problem arises naturally in frame theory and in sampling theory \cite{ABK}. More generally matrices with certain decay  far from the diagonal have been extensively studied over the last 20 years (see for instance \cite{B,J, FGL,FGL2,S1}). It has applications in various fields of analysis, such as  pseudo-differential operators \cite{Sj, G4},  
numerical analysis \cite{CS, S2, S3}, wavelet analysis \cite{J}, time-frequency analysis 
\cite{G1, G2, G3}, sampling \cite{ABK, CG, G3}), and Gabor frames \cite{BCHL,CG,Sj}.

\subsection{Left-invertibility of thin-sparse operators}

Recall that a discrete metric space $X$ is called doubling with doubling constant $D$ if for all $r>0$ and $x\in X$
$$V(x,2r)\leq DV(x,r),$$
where $V(x,r)$ denotes the cardinality of the closed ball of radius $r$. Examples of doubling metric spaces are $\Z^n$, and more generally groups with polynomial growth. 
Recall that a countable group $G$ has polynomial growth if for every finite subset $U\subset G$, there exists $C=C(U)$ and $d=d(U)$ such that $|U^n|\leq Cn^d$. By a deep theorem of Gromov \cite{Gro}, a finitely generated group $G$ has polynomial growth if and only if has a nilpotent normal subgroup of finite index. It then follows from \cite{Gui} that there exists an integer $d=d(G)$ such that for all finite symmetric generating subset $U$ of $G$, there exists $C=C(U)$ such that 
$$C^{-1}n^d\leq |U^n|\leq Cn^d.$$  
As a result, the group $G$, equipped with the {\it word metric} $d_U(g,h)=\inf\{n\in \N, g^{-1}h\in U^n\}$ is a doubling metric space.

Given a doubling metric space $X$ and a countable set $Y$, we consider an operator $A=(a_{y,x})_{(y,x)\in Y\times X}$, bounded on $\ell^2$,
whose rows are supported in balls of bounded radius (i.e. are {\it thin}), and whose columns have only a bounded number of non-zero entries (i.e. are {\it sparse}): we call such a matrix
{\it thin-sparse}.

Our first main result states that if $A$ is bounded below in $\ell^p$ for some $1\leq p\leq \infty$,
then, $B=(A^*A)^{-1}A^*$ defines a left-inverse for $A$, which is uniformly bounded on $\ell^q$ for $q\in [1,\infty]$.

\begin{thm}\label{mainThmIntro}
Let $X$ be a doubling metric space and
let $A=(a_{y,x})_{(y,x)\in Y\times X}$ be thin-sparse matrix with bounded coefficients.  Then,
\begin{itemize}
\item either
$$\lambda_p(A)=0$$ for all $1\leq p\leq \infty,$
\item or there exists $C<\infty$, such that $B=(A^*A)^{-1}A^*$ satisfies
$$\|B\|_{p\to p} \leq C,$$ for all $1\leq p\leq \infty,$
and hence defines a left-inverse for $A$.
\end{itemize}
\end{thm}

\begin{rem}
Note that for a matrix $A$ whose rows have bounded support, a
uniform bound on the coefficients is equivalent to the fact that
$A$ is bounded in $\ell^{\infty}.$ So, if $A$ is bounded in
$\ell^p$ for some $1\leq p\leq \infty$, as in particular its
coefficients are bounded, it is also bounded in $\ell^{\infty}$.
Hence by interpolation, it is bounded for all $p\leq q\leq
\infty$.
\end{rem}
We shall discuss the optimality of this result latter in subsection \ref{OptimalitySection}. 
One can actually drop the assumption of sparseness on the columns of $A$, and obtain the following stronger statement (indeed Theorem~\ref{mainThmIntro} follows by taking $p<1$ in the following theorem). Say that a matrix $(a_{y,x})_{(y,x)\in Y\times X}$ is thin-$\O$ if rows are thin, i.e. supported on balls of bounded radius (and no assumption is made on columns). 
\begin{thm}\label{mainThmIntro'}
Let $A=(a_{y,x})_{(y,x)\in Y\times X}$ be a thin-$\O$ matrix. Assume moreover that $A$ is bounded as an
operator $\ell^p(X)\to \ell^p(Y)$ for some $0< p< \infty$
(equivalently bounded on $\ell^q$ for all $p\leq q\leq \infty$). Then,
\begin{itemize}
\item either
$$\lambda_q(A):=\inf_{f\neq 0}\frac{\|Af\|_q}{\|f\|_q}=0$$ whenever $p < q\leq \infty$ and $q\geq 1;$
\item or there exists $c>0$, such that
$$\lambda_q(A)\geq c,$$ if $\max(p,1)\leq q\leq \infty.$
 In the latter case, if $p\leq 2$, then $B=(A^*A)^{-1}A^*$ defines a left-inverse for $A$, which is uniformly bounded on $\ell^q$ for $$\max(p,1)\leq q\leq p/(\max(p,1)-1).$$
\end{itemize}
\end{thm}
The conclusion of Theorem~\ref{mainThmIntro'} is optimal as one
can easily construct for every $1\leq p\leq \infty$ a matrix
$A=(a_{y,x})_{y,x\in \N}$ with one non-zero coefficient in each
row and such that
\begin{itemize}
\item $A$ is bounded in $\ell^q$, for $q\geq p$, \item
$\lambda_p(A)>0$, \item $\lambda_q(A)=0$ for all $p<q\leq \infty$.
\end{itemize}
To see this, consider a matrix such that the $n$'th column
contains exactly $n$ non-zero coefficients equal to $n^{-1/p}$,
such that the columns are piecewise orthogonal (i.e. have disjoint  supports).

\begin{rem}
Theorem~\ref{mainThmIntro} has been proved recently
\cite{ABK} for slanted matrices: let $\alpha\in \R^*$, a matrix
$(a_{y,z})_{y,z\in \Z^d}$ is called $\alpha$-slanted if its
support in $\Z^d\times\Z^d$ lies at bounded distance from the
subspace of $\R^d\times\R^d$ defined by
$\{(x,y)\in\R^d\times\R^d,\; y=\alpha z\}$.
\end{rem}

Although our proof is clearly different from the one of \cite{ABK}, both approaches share an important
idea which consists in restricting $A$ to functions supported in
balls of radius $L$. This reduces the problem to dimension
$\lesssim L^d$, which enables us to use quantitative comparisons between
$\ell^p$-norms, before letting $L$ go to infinity. Precisely, we prove the
following fact which might be of independent interest (see Theorem~\ref{MainThmApprox} for a more general
statement).
\begin{thm}\label{mainThmIntroApprox}

Let $X$ be a doubling metric space, and
let $A=(a_{y,x})_{(y,x)\in Y\times X}$ be a thin-$\O$ matrix. Assume that the matrix $|A|=(|a_{y,x}|)_{y\in Y,x\in X}$ defines a bounded operator $\ell^{p}(X)\to
\ell^{p}(Y)$, for some $1\leq p\leq \infty$. Then, there exist
$C_1$ and $C_2$ such that for all $L\geq 1$, there is a non-zero
function $h$ supported in a ball of radius $L$ such that for
all $p\leq q\leq  \infty$,
$$\frac{\|Ah\|_q}{\|h\|_q}\leq C_1\lambda_q(A)+\frac{C_2}{L}.$$
($C_1$ only depends on the space $X$, and for $X=\Z$, we can take
$C_1=6$. But $C_2$ also depends on $\||A|\|_{p\to p}$).
\end{thm}
The estimate in $O(1/L)$ for the error term is optimal as one can
easily check with $A=1-P$, where $P$ is\footnote{Note that $P$ is the diffusion operator associated with the simple random walk on $\Z$.}
 the convolution by the normalized characteristic
function of $\{-1,1\}$, acting on $\ell^p(\Z)$. 

\subsection{Application to Schur operators}

We are able (see Theorem \ref{t-thinTheorem}) to extend Theorem \ref{mainThmIntro} in a way to include all matrices which can be approximated in a suitable sense by thin-sparse matrices.  Here, we only focus on a special case, i.e. where $X=Y$ and where the matrices can be approximated by banded ones.

We will say that a matrix $(a_{x,y})$ indexed by a metric space $X$ is $N$-banded (or has propagation $\leq N$) if
$a_{x,y}=0$ as soon as $d(x,y)>N$.

We will denote by $\A$ the algebra of Schur operators.
Recall a Schur operator on $\ell^2$ is an operator which is bounded both on $\ell^1$ and on $\ell^{\infty}$, its Schur norm being defined as $\|A\|_{\A}=\|A\|_{1\to 1}+\|A\|_{\infty\to \infty}=\sup_{i}\sum_{j}|a_{i,j}|+\sup_{j}\sum_{i}|a_{i,j}|.$

\begin{thm}\label{approxThm}
Let $X$ be a doubling metric space, and let $A=(a_{x,y})$ be a Schur matrix indexed by $X$ such that there exists a sequence of $r$-banded matrices $A_r$ such that 
$$r^{t}\cdot\|A-A_r\|_{\A}\overset{r\to \infty}{\longrightarrow} 0,$$
for some $t>0$. Then the following are equivalent
\begin{itemize}
\item $A$ is bounded below for some $1\leq p\leq \infty$, 
\item $A$ is bounded below for all such $p$,
\item $B=(A^*A)^{-1}A^*$ defines a left-inverse of $A$ lying in $\A$.
\end{itemize}
\end{thm}

The first notion of weighted Schur algebra has been introduced in \cite{GL}, and then generalized in \cite{Sun}. 
Following \cite[Section 2.2]{Sun}, if $X$ is a metric space and $\omega: X\times X\to [1,\infty)$ is an admissible weight in the sense of \cite{GL}  or of \cite{Sun}, then we can define the weighted Schur algebra $\A_{\omega}(X)$ as the space of operators which are bounded for the norm
$$\|A\|_{\A,\omega}=\sup_{x}\sum_{y}\omega(x,y)|a_{x,y}|+\sup_{y}\sum_{x}\omega(x,y)|a_{x,y}|.$$
Typical admissible weights are $$\omega(x,y)=1+d(x,y)^{\alpha},$$for $\alpha\geq 0$, and $$\omega(x,y)=\exp(Cd(x,y)^{\delta}),$$ for some $C>0$, and $0<\delta<1$.  
Since the notion of admissible weight is very technical, and will never be used here, we will not recall it (or else, we suggest the reader to consider the two  previous typical examples as a definition of admissible weights since they both satisfy the conditions of \cite{GL} and of \cite{Sun}). 
\begin{cor}\label{SchurThm}
Let $X$ be a doubling metric space, and let $\omega$ be an admissible weight such that 
$\omega(x,y)\geq d(x,y)^{\alpha}$ for some $\alpha>0$. Then the following are equivalent
\begin{itemize}
\item $A$ is bounded below for some $1\leq p\leq \infty$, 
\item $A$ is bounded below for all such $p$,
\item $B=(A^*A)^{-1}A^*$ defines a left-inverse of $A$ lying in $\A_{\omega}(X)$.
\end{itemize}
\end{cor}
\bpr First an easy observation shows that the matrices $A_N$ obtained naïvely by replacing all coefficients $a_{x,y}$, where $d(x,y)>N$ by zeros satisfy the hypothesis of Theorem \ref{approxThm}. 
The last statement follows from Theorem \ref{approxThm}, together with the facts that $\A_{\omega}(X)$ is an involutive algebra, and is spectral (or inverse-closed), which are both proved in \cite{GL,Sun} (for different types of weights). Namely, since $\A_{\omega}(X)$ is involutive, $A^*\in \A_{\omega}(X)$, as it is an algebra, $A^*A\in \A_{\omega}(X)$, since it is spectral, 
$(A^*A)^{-1}\in \A_{\omega}(X)$, and finally, we conclude using that $\A_{\omega}(X)$ is an algebra.
 \epr

\subsection{Application to the class of convolution-dominated operators}

Let $G$ be a discrete group. Recall the Gohberg-Baskakov-Sjöstrand class  \cite{Sun} (also called the  convolution dominated operators class \cite{FGL2}) $\CC(G)$ is the set of all operators on $\ell^2(G)$ which are bounded for the following norm 
$$\|A\|_{\CC(G)}=\sum_{k\in G}\sup_{g^{-1}h=k}|a_{g,h}|.$$
Let $\omega$ be an admissible weight. We shall also suppose that $\omega$ is left-invariant, i.e. satisfies\footnote{Observe that the two typical classes of weights defined at the previous subsection are indeed left-invariant, when defined with a left-invariant metric.} $\omega(gk,gh)=\omega(k,h)$ for all $g,h,k\in G$. Following \cite{FGL2}, one can define the weighted convolution dominated algebra, comprising all matrices $A$ which are bounded for the following norm 
$$\|A\|_{\CC_{\omega}(G)}=\sum_{k\in G}\sup_{g^{-1}h=k}\omega(g,h)|a_{g,h}|.$$

\begin{thm}\label{CDThm}
Let $G$ be a group with polynomial growth, and let $\omega$ be an admissible left-invariant weight such that 
$\omega(g,h)\geq d(g,h)^{\alpha}$ for some $\alpha>0$. Then the following are equivalent
\begin{itemize}
\item $A$ is bounded below for some $1\leq p\leq \infty$, 
\item $A$ is bounded below for all such $p$,
\item $B=(A^*A)^{-1}A^*$ defines a left-inverse of $A$ lying in $\CC_{\omega}(G)$.
\end{itemize}
\end{thm}
The proof is completely similar to that of Theorem \ref{SchurThm} using the fact, proved in \cite{FGL2} (see also \cite{Sun} for a weaker statement) that $\CC_{\omega}(G)$ is a spectral involutive algebra for all admissible weight.

It turns out that our condition on the weight is not optimal. Indeed, in a very recent paper, Shin and Sun  managed to prove the above theorem for any admissible weight when $G=\Z^n$ \cite{Sun2}.
We believe that their proof should also work for a group with polynomial growth, although this remains to be checked carefully.

Finally, let us mention that even in the context of convolution operators on a group of polynomial growth, 
the above theorem is new, and has the following application. In view of \cite[Theorem 4.3]
{Ch}, we obtain

\begin{cor}
Let $G$ be a group with polynomial group, and suppose that an element $A\in \C G$ is bounded below in $\ell^p$ for some $1\leq p\leq \infty$, then $A$ is invertible in $B(\ell^q(G))$ for all $1\leq q\leq \infty$. \epr
\end{cor}

\subsection{Optimality of the assumptions of Theorem \ref{mainThmIntro} and Corollary \ref{SchurThm}}\label{OptimalitySection} 

There are two natural questions arising from Corollary \ref{SchurThm}. Namely, can we relax, or simply drop  one of the two main assumptions: the doubling condition on the space $X$, and the strict polynomial decay of the coefficients?

First, Corollary \ref{SchurThm} cannot be extended to the unweighted Schur algebra $\A$ since we exhibited in \cite{T} a matrix in $\A$ which is bounded below in $\ell^2$ but not in $\ell^{\infty}$. 
As Nigel Kalton pointed to me, this fact is actually well-known amongst interpolation theoretists. An easy example is $A=I-D$, where $D$ is the dilation operator on $\ell^2(\N)$, i.e. $$D(a_0,a_1,\ldots)=(a_0/2,a_0/2,a_1/2,a_1/2,\ldots).$$
Note that the operator $A^*=1-D^*$ is invertible in $\ell^2$ but not left-invertible in $\ell^1$. Indeed, the sequence of normalized characteristic functions $\phi_n=1_{[0,n-1]}/n$ satisfies $\|A^*\phi_n\|_1\to 0$.
One can extend this idea to get examples which are not left-invertible in $\ell^p$ for $1<p<2$, by replacing $D$ by $\lambda D$, where $1<\lambda<\sqrt{2}.$

Note that these examples do not exhibit any decay at infinity. On the other hand, the example given in \cite{T} is a banded matrix\footnote{Indeed, the operator considered in \cite{T} is a symmetric element of the group algebra of the free group with two generators $F_2$ seen as a convolution operator on $\ell^p(F_2).$}
indexed by the vertex set of the $3$-regular tree $T$. Therefore it belongs to $\A_{\omega}(T)$ for {\it any} weight $\omega$ on $T$. Hence, it gives a partial answer to the question of whether the metric space is required to be doubling or not. Actually, it is easy to see that $T$ has exponential growth, and therefore does not satisfy the doubling condition. Moreover, as we will see below, $T$ is a key example among those spaces\footnote{Indeed, it is an open question whether a discrete metric space $X$ with exponential growth admits a Lipschitz embedded copy of $T$.}. Note that a matrix indexed by $T$ can be easily  ``extended" to a matrix indexed by $X$ still satisfying the properties we are interested in. This provides a wide class of examples of metric spaces for which Corollary \ref{SchurThm} (and actually even Theorem \ref{mainThmIntro} for banded matrices) fails to be true. 
For instance, this excludes any metric space which is 
the vertex set of some non-amenable $k$-regular graphs. Those are graphs satisfying an isoperimetric inequality 
$$|\partial A|\geq c|A|,$$ 
for every finite subset  $A$ of vertices of the graph, where $c$ is some positive constant. The boundary $\partial A$ denotes the set of edges joining vertices of $A$ to its complement. Indeed, by the main result of \cite{BS}, such a graph admits a bi-Lipschitz embedded 3-regular tree. 
Most known finitely generated groups have exponential growth, and among them, a large class have been shown to admit a Lipschitz embedded copy of $T$: this comprises by the previously mentioned result the huge class of non-amenable groups, while for instance Rosenblatt \cite{R} proved it for non-virtually nilpotent solvable groups, which form a large class of amenable groups with exponential growth.

However, there is still an interesting question which remains open: sticking to matrices indexed by $\Z$ for instance, does the conclusion of Corollary \ref{SchurThm} hold for --say--  logarithmic decay?

\subsection{About the proofs}

The proofs of Theorem~\ref{mainThmIntro} and of its variants  split into two main parts. First, we need to show that if $A$ is bounded below for some $p$, then it is uniformly bounded below in $\ell^q$ for all $q$'s. The second part of the proof consists in showing that the left-inverse exists and is uniformly bounded in $\ell^p$ for all $p$'s. Let us now explain how the second part follows from the first one. We will deduce it from the following elementary observation. 
\begin{prop}\label{generalityprop}
Let $X$ and $Y$ be two sets, and let $A$ be an operator $\ell^2(X)\to \ell^2(Y)$ such that $A$ and $A^*$ are uniformly bounded in $\ell^p$ for all $1\leq p\leq \infty$. We have
\begin{itemize}
\item $\lambda_2(A^*A)=\lambda_2(A)^2$,
\item if $A$ is self-adjoint and $\lambda_p(A)>0$, for all $1\leq p\leq \infty$, then $A$ is invertible in $\ell^p$, and
$\|A^{-1}\|_p=1/\lambda_p(A).$
\end{itemize}
\end{prop}
\bpr The first statement simply follows from $$\lambda_2(A^*A)=\inf_{\|f\|_2=1}\langle A^*Af,f\rangle=\inf_{\|f\|_2=1}\|Af\|_2^2=\lambda_2(A)^2.$$
To show the second statement, observe that since $A$ is self-adjoint, $\lambda_2(A)>0$ implies that $A$ is invertible in $\ell^2$. Hence, $A^{-1}$ is defined on $\ell^p(Y)\cap\ell^2(Y)$ which is dense in $\ell^p(Y)$ for all $p$. But then 
\begin{eqnarray*}
\lambda_p(A)  &  =  & \inf_{f\in \ell^p(Y)\cap\ell^2(Y)}\frac{\|Af\|_p}{\|f\|_p}\\
                          & =  & \inf_{f\in \ell^p(Y)\cap\ell^2(Y)}\frac{\|f\|_p}{\|A^{-1}f\|_p}\\
                          &  =  & 1/\|A^{-1}\|_{p\to p}. 
\end{eqnarray*}
So the proposition is proved. \epr

\

To fix the ideas, let us focus on the second statement of Theorem \ref{mainThmIntro}, assuming the first statement.
If $\lambda_p(A)\geq c>0$ for all $1\leq p\leq \infty$, then in particular, this is true for $p=2$. So $\lambda_2(A^*A)\geq c^2$, which implies that $A^*A$ is invertible. But $\lambda_2(A^*A)>c^2$, and by  Proposition~\ref{MultiplicationRuleProp},  $A^*A$ is banded. So by the first statement of Theorem \ref{mainThmIntro} applied to $A^*A$, there exists $c'>0$ such that $\lambda_p(A^*A)\geq c'$ for all $1\leq p\leq \infty.$
Finally as $\|(A^*A)^{-1}\|_p=1/\lambda_p(A^*A)\leq 1/c'$, we conclude that $B=(A^*A)^{-1}A^*$ satisfies
$$\|B\|_p\leq \|A^*\|_{p}/c',$$
which is bounded independently of $p$.

\begin{rem}
Note that the fact that  the left-inverse $A^*(A^*A)^{-1}$ is uniformly bounded in $\ell^p$ for all $p$ is also an immediate consequence of the fact that $(A^*A)^{-1}$ lies in the Schur algebra \cite{GL,Sun}. 
\end{rem}

Let us now summarize  the first part of the proof of Theorems~\ref{mainThmIntro}, \ref{mainThmIntro'}. Let us assume that $\lambda_{p_0}>0$ for some $1\leq p_0\leq \infty$. In views of Proposition \ref{generalityprop}, we only need to show that $\lambda_p>0$ for all $p$.
\begin{enumerate}
\item
The first step, Theorem~\ref{mainThmIntroApprox}, is the
central part of this paper (see Section~\ref{ApproxSection}).  We show that the doubling property can be used to approximate the $\ell^p$-norm of a function $f$ by taking the norm of its projection over a subset consisting of a union of distant balls of fixed radius. However, the naive idea consisting in applying $A$ directly to this projection would only yield an error term in $L^{1/p}$, which would not enable us to deduce anything from the statement that $\lambda_{\infty}(A)>0$ (but would work for any $p<\infty$). Instead, we multiply $f$ by a certain Lipchitz function which is also supported on a union of distant balls.

\item
To obtain the uniform lower bound for $\lambda_q(A)$, using Theorem~\ref{mainThmIntroApprox}
is quite technical but the general idea is easy to understand: Theorem~\ref{mainThmIntroApprox} says that we can approximate $\lambda_q(A)$ by quotients of the form $\frac{\|Ah\|_q}{\|h\|_q}$, where $h$ are supported in balls of radius $L$ (hence, restricting to subspaces of dimension $\approx v(L)$ which is roughly less than $L^d$ for some $d$), and the error that we make is roughly in $1/L$. Comparing these quotients for different values of $q$ (and the same function $h$), we multiply our error term by $L^{d|1/p-1/q|}$. The resulting error term will therefore go to zero if $p$ and $q$ are close enough, namely if $d|1/p-1/q|<1$. Then, we just need to ``propagate" the comparison that we get between $\lambda_p(A)$ and $\lambda_q(A)$ to obtain a uniform lower bound. 
Note that similar ideas are used in \cite{ABK,Sun2}.

\item
Then, we
extend Theorem~\ref{mainThmIntro} to operators that are somehow ``polynomially approximated"
by thin-sparse operators: we call them almost thin-sparse operators (see Section~\ref{s-thinSection}). The idea of the proof is very similar to step 2 (see Lemma~\ref{mainlemma'}).

\item The proof of Theorem~\ref{mainThmIntro'}  essentially consists  in showing that a thin-$\O$ operator which is bounded in $\ell^p$, is almost thin-sparse in $\ell^q$ for all $q>p$, which is easily checked.

\end{enumerate}

\bigskip

\noindent \textbf{Acknowledgments.} I am grateful to Akram
Aldroubi, Ilia Krishtal, Qiyu Sun, Karlheinz Gröchenig for valuable discussions. I also thank Nigel Kalton for telling me his example of a matrix in $\A(\N)$ which is invertible in $\ell^2$ and not in $\ell^1$. I also thank Yemon Chu for pointing me his interesting paper \cite{Ch}, and for his remarks and corrections. 

\section{Notation for thin-sparse operators}
In all the sequel, $X$ and $Y$ are discrete metric spaces with
bounded geometry (balls of radius $r$ have less than $v(r)$
elements, for a given function $v$). However, in the definition of
thin-sparse operators, only $X$ needs a structure of metric space
($Y$ can be any set).

Let $C_c(X)$ be the space of finitely supported real-valued
functions on $X$. Let $A$ be a linear map from $C_c(X)$ to
$R^Y$. The kernel (also called the matrix) of $A$,
$(a_{y,x})_{(y,x)\in Y\times X}$ is defined by the relation
$$Af(y)=\sum_{x\in X}a_{y,x}f(x),$$
for every $f\in C_c(X)$. Conversely a matrix, i.e. a family of
reals $(a_{y,x})_{(y,x)\in Y\times X}$ defines a linear morphism
by the same formula.

The row of index $y\in Y$ of $A$ is the vector $(a_{y,x})_{x\in
X}$ of $\R^X$. The column of index $x\in X$ of $A$ is the vector
$(a_{y,x})_{y\in Y}$ of $\R^Y$. The support of $A$ is the subset
of $Y\times X$ on which $a_{y,x}\neq 0$. We define similarly the
support of a row or of a column of $A$.

\begin{nota}
If the rows of a matrix $A=(a_{y,x})_{y\in Y}$ satisfy some
property ``P'', and if its columns satisfy some property ``Q'', we will say
that ``$A$ is P-Q''. If we make no assumption on the columns, we will
say that $A$ is P-$\O$, and so on. We will consider two properties
for the rows or the columns:
\begin{itemize}

\item We say that the rows (or the column) of $A$ are thin, of
thickness at most $r$ if their support are contained in balls of
radius $r$.

\item  We say that the rows (or the columns) are sparse, of
sparseness at most $v$ if their support has cardinality at most $v$.

\item We denote by $TS(X,Y)$ (resp. $ST(X,Y)$, $T(X,Y)$, $\O
T(X,Y)$ and $T\O(X,Y)$) the space of thin-sparse (resp.
sparse-thin, thin-thin, $\O$-thin and thin-$\O$) operators.
\end{itemize}
\end{nota}

As the spaces have bounded geometry, sparse is a weaker condition
than thin. Hence sparse-sparse is weaker than thin-sparse, which
is weaker than thin-thin, etc.

\begin{rem}
A particular case of thin-thin matrices (when $X=Y$) are matrices for which the
support is contained in $\{(y,x)\in X^2, d(x,y)\leq r\}$ for some
$r>0$. Such matrices are sometimes called banded, or with finite propagation.
\end{rem}

\begin{nota}

\

\begin{itemize}

\item For all $1\leq p\leq \infty$, the norm of an operator
$A:\ell^p(X)\to\ell^p(Y)$ is called the $\ell^p$-norm of $A$ and
is denoted by $\|A\|_{p\to p}.$

\item Let $A=(a_{y,x})_{(y,x)\in Y\times X}$. The absolute value
of $A$ is operator $|A|=(|a_{y,x}|)_{(y,x)\in Y\times X}$.

\item We say that $A$ is absolutely uniformly bounded if
$$\sup_{1\leq p\leq \infty}\||A|\|_{p\to p}<\infty.$$
\end{itemize}
\end{nota}

\section{Preliminary remarks about thin-sparse operators}

\subsection{Combinatorial properties}

The following easy fact is a crucial property of TS operators. We say that two subsets $U$ and $V$ of a metric space are $t$-disjoint if $d(x,y)> t$ for all $(x,y)\in U\times V$.
\begin{prop}\label{support_TS_prop}
Let $X$ be a metric space, and $Y$ be a set. Let $A$ be
a thin-$\O$ operator of thickness $r$ and let $v$ and $u$ be
two functions on $X$ whose supports are $2r$-disjoint. Then, $Au$
and $Av$ (which are well defined functions) have disjoint support.
\end{prop}
\bpr We just have to consider a row $L$ of $A$ and to prove that
$\langle L,u \rangle \neq 0$ implies $\langle L,v \rangle=0$. But this is a trivial consequence of
the fact that $L$ is supported in a ball of radius $r$, which has diameter $\leq 2r$, and that
the supports of $u$ and $v$ are at distance $>2r$.\epr

\

The following proposition is straightforward and left as an
exercise.

\begin{prop}\label{MultiplicationRuleProp}
Let $X$ be a metric space and let $Y$ be a set. If $A\in T\O(X,Y)$
then $A^*A$ (when it exists) is banded.
\epr
\end{prop}

\subsection{Norms of sparse-sparse operators are equivalent}

\begin{prop}\label{normSS}
A sparse-sparse operator $A$ is absolutely uniformly bounded, if
and only if it is bounded in $\ell^p$ for some $1\leq p\leq
\infty$, if and only if it has bounded coefficients.
\end{prop}

\bpr Let $X$ and $Y$ be two sets and let $A=(a_{y,x})_{(x,y)\in
X\times Y}$ be a sparse-sparse operator of sparseness $v$. Note
that the norm $\|A\|_{\infty}=\sup_{(y,x)\in Y\times X}|a(y,x)|$
is trivially less than all operator norms. Hence it is enough to
prove that for every $1\leq p\leq \infty$, $\|A\|_{p\to p}\leq
C\|A\|_{\infty}$ for some $C$ depending only on $v$. Fix $y\in Y$,
and let $S_y$ be the support of the corresponding row
$(a_{y,x})_{x\in X}$. For every $f\in C_c(X)$,

$$|Af(y)|=\left|\sum_{x\in X}a_{y,x}f(x)\right|\leq
\|A\|_{\infty}\sum_{x\in S_y}|f(x)|.$$ Hence, using Hölder's
inequality and the majoration $|S_y|\leq v$ for all $y\in Y$, we
obtain
\begin{eqnarray*}
\|Af\|_p^p &\leq & \|A\|_{\infty}^p\sum_{y\in
Y}\left(\sum_{x\in S_y}|f(x)|\right)^p\\
& \leq & \|A\|_{\infty}^p\sum_{y\in Y}v^{p-1}\sum_{x\in
S_y}|f(x)|^p
\end{eqnarray*}
Now, note that for every $x\in X$ and every $k\in \N$, $f(x)$
appears $k$ times in the sum above if there are $k$ distinct
elements of $Y$, $y_1,\ldots, y_k$ such that $x\in
S_{y_1}\cap\ldots\cap S_{y_k}$, hence if $y_1,\ldots, y_k$ lie in
the support of the column $(a_{y,x})_{y\in Y}$. But as the
sparseness of $A$ is at most $v$, this implies that $k \leq v$.
Therefore, we have
\begin{eqnarray*}
\|A\|_{\infty}^p\sum_{y\in Y}v^{p-1}\sum_{x\in
S_y}|f(x)|^p& \leq & v^p\|A\|_{\infty}^p\sum_{x\in X}|f(x)|^p\\
& = & v^p\|A\|_{\infty}^p\|f\|_p^p.\quad \blacksquare
\end{eqnarray*}

\section{Proof of the approximation property}\label{ApproxSection}

Recall that a discrete metric space $X$
is said to be doubling of doubling constant $C<\infty$ if for all $x\in X$ and every $r>0$,
$$|B(x,2r)|\leq C|B(x,r)|.$$
Our purpose in this section is to prove the following theorem

\begin{thm}\label{MainThmApprox}
Assume that $X$ is a doubling metric space and let $A\in T\O(X,Y)$
of thickness $r$, such that $\||A|\|_{p\to p}\leq 1$ for some
$1\leq p< \infty$. There exists $C$ such that for every $f\in
L^p(X)$, and every $L\geq r$, there exists a function $h\in
L^p(X)$ supported in a ball of radius $2L$ such that
$$\frac{\|Ah\|_p}{\|h\|_p}\leq C\left(\frac{\|Af\|_p}{\|f\|_p}+\frac{r}{L}\right),$$
where, the quantity $C$ only depends on the doubling constant of $X$.
\end{thm}

\subsection{Coloring of a family of balls}

Recall that a $d$-coloring of a set $\mathcal{P}$ of subsets of $X$
is a map $$j:\mathcal{P}\to \{1,2, \ldots,d+1\}$$ such that every
two elements in $\mathcal{P}$ with the same color (i.e. same
image
by $j$) are disjoint.

Also classical is the notion of coloring of a graph: a
$d$-coloring of a graph $\mathcal{G}$ is a map
$$j:V(\mathcal{G})\to \{1,2, \ldots,d+1\},$$ where $V(\mathcal{G})$ is
the vertex set of $\mathcal{G}$, such that any two adjacent
vertices have distinct colors. A classical result of graph theory,
known as Brooks' theorem says that any graph of degree at most $d$
admits a $d$-coloring.

It tuns out that these two definitions of coloring are related via the
notion of dual graph. Recall that the dual graph $\mathcal{G}$ of
$\mathcal{P}$ is defined as follows: the set of vertices
$V(\mathcal{G})$ is $\mathcal{P}$, and two vertices are adjacent
if and only if they have a non-empty intersection. Clearly, a
$d$-coloring of $\mathcal{G}$ yields a $d$-coloring of
$\mathcal{P}$ and conversely.

We will need the following lemma. 
\begin{lem}\label{coloringlem}
Let $X$ be a doubling metric space and let $\alpha\geq 1$. There
exists an integer $d$ such that for every $L>0$, there exists a
covering of $X$ by balls of radius $L$ admitting a $d$-coloring
such that the centers of two balls of same color are at distance
$\geq \alpha L$ from one another.
\end{lem}
\bpr Consider a minimal covering $\mathcal{B}=(B(x_i,L))_i$ of $X$ (which exists since $X$ is doubling). By minimality, the
balls $B(x_i,L/4)$ are piecewise disjoint.

Now, consider the covering $\mathcal{B}'=(B(x_i,\alpha L))_i$. It
is  easy to see that the doubling property implies
that the dual graph of $\mathcal{B}'$ has degree less than a
certain constant $d$. Indeed, for every $i$, let $d_i$ be degree at the vertex $i$ of the dual graph. In other words, $d_i$ is the number of balls $B(x_j,\alpha L)$ with $j\neq i$, intersecting $B(x_i,\alpha L)$. Let $J_i$ be the set of such indices. Note that the disjoint union $\cup_{j\in J_i}B(x_j,L/4)$ is contained in $B(x_i,4\alpha L)$. On the other hand, by the doubling property, there exits $c>0$ only depending on $\alpha$ such that  $\inf_{j\in J_i}V(x_j,L/4)/V(x_i,4\alpha L)\geq c.$ But since
$$d_i\inf_{j\in J_i}V(x_j,L/4)\leq V(x_i,4\alpha L),$$
we deduce that $d_i\leq 1/c$, so that we can set $d=[1/c].$

Hence, by Brooks' theorem, this graph admits
a $d$-coloring, which means that $\mathcal{B}'$
has a $d$-coloring. Inducing this coloring to $\mathcal{B}$ yields
the desired $d$-coloring. \epr

\subsection{Approximating a function by a function supported by a disjoint union of balls of fixed radius.}

In the following lemma we characterize the doubling condition in
terms of approximation of functions by functions supported by
disjoint unions of balls of fixed radius.

For every subset $\Omega$ of a metric space $X$ and every $L>0$, we
denote
$$[\Omega]_L=\{x\in X,d(x,\Omega)\leq L\}.$$ We also denote the
characteristic function of a subset $\Omega$ by $1_{\Omega}$.
Finally, a $K$-separated subset of $X$ is a subset whose elements are pairwise at distance at least $K$.

\begin{lem}\label{approxlem}
A metric space $X$ is doubling if and only if for every
$\alpha\geq 1$, there exists a constant $c>0$ such that for every
$1\leq p\leq \infty$, every $f\in \ell^p(X)$ and every $L>0$, one
can find an $\alpha L$-separated subset $P$ of $X$ such that
$$\|1_{[P]_L}f\|_p\geq c\|f\|_p.$$
\end{lem}
\bpr Consider the covering $\mathcal{B}$ of the previous lemma and
for every $1\leq k\leq d+1$, let $P_k$ be the set of centers of
balls of $\mathcal{B}$ with same color $k$. Since
$X=\bigcup_{k=1}^{d+1}[P_k]_L$, we have
$$\|f\|_p\leq \|\sum_{k}1_{[P_k]_L}|f|\|_p\leq
\sum_{k}\|1_{[P_k]_L}f\|_p\leq (d+1)\max_k \|1_{[P_k]_L}f\|_p.$$
So Lemma~\ref{approxlem} follows taking $P=P_k$ with a $k$ for
which the max is attained. The converse follows by taking $f$ to
be the characteristic function of a ball of radius $2L$ and
$\alpha\geq 6$, so that the intersection between $[P]_L$ and our
ball of radius $2L$ is contained in a single ball of radius $L$.
\epr

\medskip

\noindent{\bf In the sequel, we fix $\alpha=6$.}

The following lemma is trivial and left to the reader.
\begin{lem}\label{function_delta_lemma}
For each $P$ like in the previous lemma, the function $\Delta_P$, defined by
$$\Delta_P(x)=\max\{0,1-d(x,P)/(2L)\},$$ satisfies
\begin{enumerate}
\item $\Delta_P=0$ outside of $[P]_{2L}$

\item $\Delta_P\geq 1/2$ on $[P]_L$.

\item $\Delta_P$ is $1/(2L)$-Lipschitz.

\item $0\leq \Delta_P\leq 1$.
\end{enumerate}\epr
\end{lem}

\begin{rem} Keeping the notation of the previous lemmas, the function
$g=\Delta_P f$ satisfies, thanks to the second property of
$\Delta_P$ and to Lemma \ref{approxlem},
$$\|g\|_p\geq c\|f\|_p.$$
On the other hand, the support of $g$ is contained in a union of
$4L$-disjoint balls of radius $2L$. Write $g=\sum_ig_i$, where
each $g_i$ is supported in one of those balls. Assume that $4L\geq
2r$. Then by Proposition~\ref{support_TS_prop},
$$\|Ag\|_p^p=\sum_i\|Ag_i\|_p^p.$$
So we have
$$\inf_i\frac{\|Ag_i\|_p}{\|g_i\|_p}\leq \frac{\|Ag\|_p}{\|g\|_p}.$$
\end{rem}

\noindent{\bf Proof of Theorem~\ref{MainThmApprox}. } Thanks to
the previous remark, we just need to prove a weaker version of the  theorem where in the conclusion, the function
$h$ is replaced by a function $g$ supported in a union of $2r$-disjoint balls
of radius $2L$. We consider $g=\Delta_Pf$, which has this property
since $L\geq r$. Let us start with a pointwise estimate. Fix some $y_0\in
Y$. For every $x,z\in X$,
$$g(x)=\Delta_P(x)f(x)=\Delta_P(z)f(x)+(\Delta_P(x)-\Delta_P(z))f(x).$$
We now specify $z=x_0$, such that the support of the row $(a_{y_0,x})_x$ is
contained in $B(x_0,r)$. We have
$$Ag(y_0)=\Delta_P(x_0)\sum_xa_{y_0,x}f(x)+ \sum_{x}a_{y_0,x}(\Delta_P(x)-\Delta_P(x_0))f(x).$$
So by Property (4) of $\Delta_P$,
$$|Ag(y_0)|\leq |Af(y_0)|+ \sum_{x}|a_{y_0,x}||\Delta_P(x)-\Delta_P(x_0)||f(x)|.$$
By Property (3) of $\Delta_P$,
$$|Ag(y_0)|\leq |Af(y_0)|+ \frac{r|A||f|(y_0)}{L}.$$
Now, taking the $\ell^p$ norm and applying the triangular inequality, we obtain
$$\|Ag\|_p\leq \|Af\|_p+ \frac{r\||A||f|\|_p}{L}.$$
We finally divide by $\|g\|_p$, and conclude thanks to the inequality $\|g\|_p\geq c\|f\|_p$.\epr

\begin{rem}\label{Zrem1} For $X=\Z$, we have $v(k)=2k+1$, and the doubling constant is less than $2$. Note that we can take $P=\{x_0+6kL,\; k\in \Z\}$ for some $x_0$. Moreover, one checks easily that a good choice of $x_0$ gives
$$\|1_{P}f\|_p\geq \|f\|_p/3.$$ Now assume that $A$ is thin-thin of thickness $\leq r$. By the proof of Proposition~\ref{normSS}, we have $\||A|\|\leq v(r)\|A\|_{\infty}$. Hence, we obtain that there exists a function $h$ supported on a ball of radius $r$ such that
$$\frac{\|Ah\|_p}{\|h\|_p}\leq 3\left(\frac{\|Af\|_p}{\|f\|_p}+\frac{3r^2\|A\|_{\infty}}{L}\right).$$
\end{rem}

\section{$\ell^p$-stability of thin-sparse operators}\label{morepreciseSection}

Here is a more general version of Theorems~\ref{mainThmIntro}, with some
precisions that we omitted in the introduction.

\begin{thm}\label{MainThm}
Let $X$ be a metric space of doubling constant $D<\infty$ and let
$Y$ be any set. Fix some $r,v>0.$ Let $A\in TS(X,Y)$ be of
thickness at most $r$, sparseness at most $v$. Assume moreover
that $\||A|\|_{p\to p}\leq 1$ for all $1\leq p\leq
\infty$. Then
there exist $c=c(r,v,D)>0$ and
$\delta=\delta(D)>0$ for all $1\leq p,q\leq \infty,$
$$\lambda_p(A)\geq c\lambda_q(A)^{\delta}.$$
\end{thm}

In Section~\ref{s-thinSection}, we prove that the conclusion
Theorem~\ref{MainThm} 
is true for more
general operators which are ``well" approximated by
thin-sparse and thin-$\O$ operators respectively.

Theorem~\ref{MainThm} (and the remark following
Theorem~\ref{MainThm'}) result from the following more precise
results. Let
$$\lambda=\inf_{p_0\leq p\leq \infty}\lambda_p(A),$$ and let 
$p_m$ be such that $\lambda_{p_m}\geq \lambda/2$.
Let 
$$\Lambda=\sup_{p_0\leq p\leq \infty}\lambda_p(A),$$ and let $p_M$ be such that $\lambda_{p_M}\leq 2\Lambda$.

Note that since $X$ is doubling, there exists $d$ and $K$ such that $V(x,R)\leq KR^d$ for all $x\in X$ and $R>0$. 

\begin{thm}\label{mainthm}
Let $A\in TS(X,Y)$ of thickness $r$, sparseness $v$ and such that
$\||A|\|_{p\to p}\leq 1$ for all $p_0\leq p\leq \infty$. Then
there exists $k=k(v,r,d)>0$ such that
$$\lambda\geq k\Lambda^{4d}.$$
\end{thm}

\begin{thm}\label{mainthmbis}
Let $A\in T\O(X,Y)$ of thickness $r$ and such that $\||A|\|_{p\to
p}\leq 1$ for all $p_0\leq p\leq \infty$. Then there exists
$k=k(r,d)>0$ such that for all $p_0\leq p\leq q\leq \infty$,
$$\lambda_p\geq k\lambda_q^{4d}.$$
\end{thm}

These theorems will be proved after a series of lemmas.

\begin{lem}\label{lambdalem}
Fix some $1\leq p_0<\infty$. Let $A\in T\O(X,Y)$ of thickness $r$
and such that $\||A|\|_{p\to p}\leq 1$ for all $p_0\leq p\leq
\infty$.
\begin{itemize}

\item[(i)] there exist $d>0$ and $C'$ (depending on the doubling
constant) such that for all $p_0\leq p\leq q\leq \infty$ and all
$L\geq r$
$$\lambda_q(A)\leq C'L^{\left|\frac{d}{p}-\frac{d}{q}\right|}\left(\lambda_p(A)+
r/L\right).$$

\item[(ii)] if moreover, $A\in TS(X,Y)$ of sparseness $v$, then
for all $p_0\leq q\leq p \leq \infty,$
$$\lambda_q(A)\leq C'v^{\left|\frac{1}{p}-\frac{1}{q}\right|}L^{\left|\frac{d}{p}-\frac{d}{q}\right|}\left(\lambda_p(A)+
\frac{r}{L}\right).$$
\end{itemize}
\end{lem}
\bpr Theorem~\ref{MainThmApprox} implies
$$\inf_{\supp(h)\subset B(x,2L)}\frac{\|Ah\|_p}{\|h\|_p}\leq C\left(\lambda_p(A)+\frac{r}{L}\right).$$
On the other hand, if $h$ is supported in a subset of size $N$,
then for $p\leq q$,
\begin{equation}\label{EquationDimension}
\|h\|_q\leq \|h\|_p\leq
N^{\left|\frac{1}{p}-\frac{1}{q}\right|}\|h\|_q.
\end{equation}
The power in $L$ appearing in the inequalities now comes from the inequality
$V(x,L)\leq KL^d$.
Indeed, if $p\leq q$, then we obtain (i) applying the left inequality of
(\ref{EquationDimension}) to $Ah$ (where the support of $Ah$ does
not play any role) and the right inequality to $h$, whose support
has cardinality at most $KL^d$. So take $C'=CK$.

If $p\geq q$, then we apply the right inequality of
(\ref{EquationDimension}) to $Ah$ for which we control the
support thanks to the sparseness of $A$'s columns. Namely, the
cardinality of the support of $Ah$ is at most $v$ times the cardinality
of $h$'s support. This explains the corresponding power of $v$ in (ii).\epr

\begin{lem}\label{lastlem}
Let $A\in T\O(X,Y)$ of thickness $r$ and such that $\||A|\|_{p\to
p}\leq 1$ for all $p_0\leq p\leq \infty$.
\begin{itemize}
\item[(i)] for all $p_0\leq p\leq \infty$, $\lambda_p(A)=0$
implies $\lambda_q(A)=0$ for all $q\geq p$. \item[(ii)] Let $K$ be twice the constant $C'$ of
Lemma~\ref{lambdalem}. Then, for all $p_0\leq p\leq q\leq
\infty,$
$$\lambda_q(A)\leq Kr^{\left|\frac{d}{p}-\frac{d}{q}\right|}\lambda_p(A)^{1-\left|\frac{d}{p}-\frac{d}{q}\right|}.$$

\end{itemize}
\end{lem}

\begin{lem}\label{lastlembis}
Let $A\in TS(X,Y)$ of
sparseness $v$ and thickness $r$, and such that $\||A|\|_{p\to
p}\leq 1$. Then, For all $p_0\leq p\leq \infty$,
\begin{itemize}
\item[(i)] For every $p_0\leq p,q\leq \infty$, $\lambda_p(A)=0$ if
and only $\lambda_q(A)=0$. \item[(ii)] Let $K$ be twice the constant $C'$ of
Lemma~\ref{lambdalem}. For all $p_0\leq p,q\leq \infty,$
$$\lambda_q(A)\leq Kv^{\left|\frac{1}{p}-\frac{1}{q}\right|}r^{\left|\frac{d}{p}-\frac{d}{q}\right|}\lambda_p(A)^{1-\left|\frac{d}{p}-\frac{d}{q}\right|}.$$
\end{itemize}
\end{lem}

\bpr Both lemmas are proved in the same way: so let us show
Lemmas~\ref{lastlembis}. To obtain (ii), take $L=r/\lambda_p(A)$
in Lemma~\ref{lambdalem}. To prove (i), we just have to note that
the vanishing of $\lambda_p(A)$ ``propagates" thanks to
Lemma~\ref{lambdalem}: $\lambda_p(A)=0\Rightarrow \lambda_q(A)=0$
if $|\frac{d}{q}-\frac{d}{p}|\leq 1/2$ (let $L\to \infty$). \epr

\

\bpr To show Theorems~\ref{mainthm} and~\ref{mainthmbis}. we
``propagate" the inequalities (ii) of Lemmas~\ref{lastlem} and
\ref{lastlembis}. As the proofs are the same for both theorems,
let us focus on the first one. If
$\left|\frac{d}{p}-\frac{d}{q}\right|\leq 1/2$, the inequality
(ii) of Lemma~\ref{lastlembis} yields
$$\lambda_p(A)\leq C(v,r,d)\lambda_q(A)^{2}.$$
Now, as $\left|\frac{d}{p_m}-\frac{d}{p_M}\right|\leq d$, we just
need to iterate this $2d$ times, which gives the theorem. \epr

\begin{rem}\label{Zrem2}
Here, assume that $X=Y=\Z$, and that $A$ is thin-thin of thickness $r$. Instead of assuming that $\||A|\|=1$, we prefer to write Lemma~\ref{lastlembis} with respect to $\|A\|_{\infty}$ (which is easier to compute in general): a consequence is that we have to replace $r$ by $3r^3\|A\|_{\infty}$.
From Remark~\ref{Zrem1} that we can take $C'=9$ in Lemma~\ref{lambdalem} (as $v(r)\leq 3r)$. Hence we can take $K=18$.
Directly from Lemma~\ref{lastlembis} (ii), we obtain that
$$\lambda_2(A)\geq \frac{\Lambda^2}{162r^{3}\|A\|_{\infty}}.$$
\end{rem}

\section{Extension to $(t,s)$-almost thin-sparse operators}\label{s-thinSection}

\begin{defn}
Fix some $t,s>0$ and some $1\leq p\leq \infty$. An operator is
$(t,s)$-almost thin-sparse for in $\ell^q$ for all $q\geq p$ if
there exists $K<\infty$ such that for all $r,v>0$, there is an
element $A_{r,v}\in TS(X,Y)$ of thickness $\leq r$ and sparseness
$\leq v$ such that $\||A-A_{r,v}|\|_{q\to q}\leq K(r^{-t}+v^{-s})$
for all $q\geq p$.
\end{defn}

This section is devoted to the proof of the following result.

\begin{thm}\label{t-thinTheorem}
Fix some $t,s>0$ and some $1\leq p_0\leq \infty$. Let $X$ be a
metric space with the doubling property, and let $Y$ be any set. Let $A$
be $(t,s)$-almost thin-sparse in $\ell^p$ for all $p\geq p_0$.
Then either $\lambda_p(A)=0$ for all $1\leq p_0\leq p \leq
\infty$, or there exists $c>0$ such that
$\lambda_p(A)>c$ for all $p_0\leq p\leq \infty$.
\end{thm}

This will result from the following analogue of Theorem~\ref{mainthm} for
$(t,s)$-almost thin-sparse operator. 
Theorem \ref{Main(t,s)-thinThm} will also be used in the proof of Theorem~\ref{MainThm'}.

\begin{thm}\label{Main(t,s)-thinThm}
Fix some $t,s>0$ and some $1\leq p_0\leq \infty$. Let $X$ be a
doubling metric space with doubling constant $D$, and let $Y$
be any set. Let $A$ be $(t,s)$-almost thin-sparse in $\ell^p$ for
all $p\geq p_0$. Then there is $c=c(D,t,s)>0$, and
$\delta=\delta(D,t,s))>0$ such that for all $p_0\leq p,q\leq
\infty,$
$$\lambda_p(A)\geq c\lambda_q(A)^{\delta}.$$
\end{thm}

In the sequel, $a\lesssim b$ will mean $a\leq Cb$, where $C=C(D,t,s)$.

\bpr First, we need the analogue of
Theorem~\ref{MainThmApprox}.

\begin{lem}\label{mainlemma'}
For all $p_0\leq p\leq
\infty$,   $f\in L^p(X)$, all $L\geq 1$ and 
$r,v>0$, there exists a
function $h\in L^p(X)$ supported in a ball of radius $2L$ such
that
$$\frac{\|Ah\|_p}{\|h\|_p}\lesssim \frac{\|Af\|_p}{\|f\|_p}+\frac{r}{L}+r^{-t}+v^{-s}.$$
\end{lem}

\bpr This is immediate, writing $A=A_{r,v}+(A-A_{r,v})$ where
$A_{r,v}$ is thin-sparse of thickness $r$ and sparseness $v$, and
using $\||A-A_{r,v}|\|_p\leq K(r^{-t}+v^{-s})$. \epr

Then we need the analogues of Lemma \ref{lambdalem}, \ref{lastlem} and \ref{lastlembis}.

\begin{lem}\label{ref'}
For all $p_0\leq p, q\leq \infty,$ and
for all $L\geq 1$ and 
$r,v>0$,
$$\lambda_q(A)\lesssim v^{|\frac{1}{p}-\frac{1}{q}|}L^{|\frac{d}{p}-\frac{d}{q}|}\left(\lambda_p(A)+\frac{r}{L}
+r^{-t}+v^{-s}\right).$$
\end{lem}
\bpr This is proved exactly as we proved
Lemma~\ref{lambdalem}.\epr

\begin{lem}\label{lastlembibis}
There exists  $u=u(D,s,t)$ such that for all $p_0\leq p, q\leq \infty,$   
$$\lambda_q(A)\lesssim \lambda_p(A)^{1-|\frac{2d}{up}-\frac{2d}{uq}|}.$$
\end{lem}
\bpr 
The proof follows by choosing in the previous lemma, $r=L^{1/2}$, $v=L^d$, and $L=\lambda_p^{-1/u}$, where $u=\min\{1/2,t/2,sd\}.$ \epr

\medskip

The proof of Theorem \ref{Main(t,s)-thinThm} now relies on an argument of propagation similar to the one used in the proof of Theorem~\ref{mainthm}. \epr

\section{Left-invertibility of thin-$\O$-operators}

\begin{thm}\label{MainThm'}
Let $X$ be a metric space of doubling constant $D<\infty$ and let
$Y$ be any set. Let $A=(a_{y,x})_{(y,x)\in Y\times X}$ be a
thin-$\O$ matrix. Assume moreover that $A$ is bounded as an
operator $\ell^{p_0}(X)\to \ell^{p_0}(Y)$ for some $0< p_0\leq
\infty$.  Then for every $p_1>p_0$, there exists
$c=c(p_1-p_0,r,D)>0$ and $\delta=\delta(p_1-p_0,D)>0$ such that
for all $\max\{1,p_1\}\leq p,q\leq \infty,$
$$\lambda_p(A)\geq c\lambda_q(A)^{\delta}.$$
\end{thm}
\begin{rem}
Before proving the theorem,  we point out that one cannot improve the theorem to have $p_0=p_1$. 
Indeed, in the spirit of the example explained in the introduction, for $r=1$ and $X=Y=\Z$, we can find a sequence of
thin-sparse operators $A_n=(a_{y,x})_{(y,x)\in Y\times X}$ of
thickness $1$, sparseness $n$, and such that
\begin{itemize}
\item $\|A_n\|_{p_0\to p_0}=\lambda_{p_0}(A_n)=1$ for all $n\in
\N$, \item and $\lambda_p(A_n)\to 0$ when $n\to 0$ for all
$p>p_0.$
\end{itemize}
On the other hand, it is interesting to note that (in virtue of Theorem \ref{mainthmbis}) there exists $c'=c(r,D)>0$ and $\delta'=\delta'(D)>0$ such
that for all $p_0\leq p\leq q\leq \infty$
$$\lambda_p(A)\geq c'\lambda_q(A)^{\delta'}.$$
\end{rem}

Theorem~\ref{MainThm'} results from
Theorem~\ref{Main(t,s)-thinThm} and from the fact that thin-$\O$
operators that are bounded in $\ell^p$ are $(1,1/p-1/q)$-almost
thin-sparse in $\ell^q$ for all $q>p$. This is a consequence of the following
proposition.

\begin{prop}
Let $X=(X,d)$ be a metric space such that balls of radius $r$ have
cardinality at most $v(r)$, and let $Y$ be a set. Fix some $\eps>0$
and some $r\geq 1$. Let $A=(a_{y,x})_{(y,x)\in Y\times X}$ be a
thin-$\O$ operator of thickness $\leq r$ such that $\|A\|_{p\to
p}=1$ for some $0< p<\infty$. Then, there is $C=C(\eps)$ such
that for every $q\geq p+\eps$ and every $m\in \N$, there exists a
thin-sparse operator $A_m$ of thickness $\leq r$, sparseness $\leq
m$ such that $$\||A-A_m|\|_{q\to q}\leq
\frac{Cv(r)^{1-1/q}}{m^{1/p-1/q}}.$$
\end{prop}
\bpr First, let us prove the following lemma.

\begin{lem}\label{lemOrderedSequences}
Let $n$ be a positive integer, and $0<a_n\leq \ldots\leq a_1$ such
that $\sum_{i=1}^na_i^p=1$, then for all $0\leq m\leq n$, and $q\geq p,$
\begin{equation}\label{eq1}
\left(\sum_{i=m+1}^{n}a_i^q\right)^{1/q}\leq
\frac{(p/q)^{1/q}(1-p/q)^{1/p-1/q}}{m^{1/p-1/q}}.
\end{equation}
In particular, for every $\eps>0$ there exists $C=C(\eps)$ such
that for all $q\geq p+\eps$,
$$\left(\sum_{i=m+1}^{n}a_i^q\right)^{1/q}\leq \frac{C}{m^{1/p-1/q}}.$$
\end{lem}
\noindent{\bf Proof of the lemma.} Let us find the maximum of the
function
$$\theta_{m,q}(a_1,\ldots, a_n)=\sum_{m+1}^n a_i^q,$$
under the conditions
$$\sum_{i=1}^n a_i^p=1,$$ and for all $1\leq i\leq n-1,$
$$a_{i+1}-a_i\leq 0.$$
\begin{clai}
The maximum of $\theta_{m,q}$ is attained at $(a_1,\ldots,a_n)$
such that $a_i=0$ for $i\geq k$ and $a_i=1/k^{1/p}$ for $i<k$, where
$k$ is an integer $\geq m+1$.
\end{clai}
\noindent{\bf Proof of the claim.} First, note that since $(a_i)$ is non-increasing, the maximum will be attained
when $a_i=a_j$ for all $i\leq j\leq m.$

On the other hand, a straightforward application of Lagrange
multipliers shows that $\theta_{m,q}$ cannot reach its maximum at a point $(a_1,\ldots,a_n)$ such that  $0<a_{i+1}<a_{i}$ for some $1\leq i \leq n-1$. Hence, if $a_{i+1}<a_i$, then $a_{i+1}=0$. There exists therefore only one such $i$. Let $k:=i+1$. Note that $\theta_{m,q}$ is not identically zero: hence, since the sequence $(a_j)$ corresponds to a maximum of $\theta_{m,q}$, $k$ has to be $\geq m+1$. Summarizing this discussion, there exists $k\geq m+1$ such that the sequence
$a_i=0$ for $i\geq k$ and $a_i=1/k^{1/p}$ for $i<k.$ 
\epr

With the notation of the claim, we have
\begin{equation}\label{eq2}
\max \theta_{m,q}=\frac{k-m}{k^{q/p}}.
\end{equation}
To finish the proof of the Lemma, note that the derivative of
$\frac{k-m}{k^{q/p}}$ with respect to $k$ vanishes exactly at the
value $m/(1-p/q)$, which corresponds to a maximum. Replacing $k$
by this value in (\ref{eq2}) yields (\ref{eq1}). \epr

\

Now, let us prove the proposition. As $\|A\|_{p\to p}=1$, for
every $x\in X$, the column $C_x=(a_{y,x})_{y\in Y}$ has
$\ell^p$-norm at most $1$. By Lemma~\ref{lemOrderedSequences},
there exists a subset $S_x$ of $Y$ of cardinality $\leq m$ such that
$$\sum_{y\in Y\smallsetminus S_x}|a_{y,x}|^q\leq C^q/m^{q/p-1}.$$
Now, we define $A_m$ from $A$ by replacing the coefficient
$a_{y,x}$ by $0$ whenever $y\in Y\smallsetminus S_x$. By
construction, $A_m$ is thin-sparse of thickness $\leq r$ and
sparseness $\leq m$.

Let $f\in \ell^q(X)$. Denote by $C_m=|A-A_m|=(c_{y,x})_{(y,x)\in
Y\times X}.$ Using Hölder inequality (which is possible since $q\geq 1$), we obtain
\begin{eqnarray*}
\||A-A_m|f\|_q^q &= &\sum_{y\in Y}\left(\sum_{x\in
X}c_{y,x}f(x)\right)^q\\
                 & \leq & \sum_{y\in Y}v(r)^{q-1}\left(\sum_{x\in
                 X}c_{y,x}^q|f(x)|^q\right)\\
                 & = & v(r)^{q-1}\sum_{x\in X}|f(x)|^q\sum_{y\in Y\smallsetminus S_x}|a_{y,x}|^q \\
                 & \leq & \frac{C^qv(r)^{q-1}}{m^{q/p-1}}\|f\|_q^q. \;
                 \blacksquare
\end{eqnarray*}

\bigskip
\footnotesize


\begin{thebibliography}{KM98b}
\bibitem[ABK]{ABK} A. {\sc Aldroubi}, A. {\sc Baskarov}, I. {\sc Krishtal}. \newblock {\em Slanted matrices, Banach frames and sampling.} \newblock  J. Funct. Anal., 255, 1667-1691, 2008.


\bibitem[Bar]{Barn}  B. A. {\sc Barnes}. \newblock {\em When is the spectrum of a convolution operator of $L\sp p$ independent of $p$?} \newblock Proc. Edinburgh Math. Soc. (2) 33, no. 2, 327-332, 1990.


\bibitem[Bas]{B} A. {\sc Baskarov}. \newblock {\em Asymptotic estimates for elements of matrices of inverse operators.} \newblock Siberian. Math. J. 38 (1), 10-22, 1997.

\bibitem[BCHL]{BCHL} R. {\sc Balan}, P.G. . {\sc Cassazza}, C. {\sc Heil} and Z. {\sc Landau}. {\em Density, overcompleteness and localization 
of frames I.} \newblock Theory; II. Gabor system, Preprint, 2004 


\bibitem[BS]{BS} I. {\sc Benjamini},  O. {\sc Schramm}, {\em Every graph with a positive Cheeger 
constant contains a tree with a positive Cheeger constant}. Geom. Funct. 
Ann. 7, 403-419, 1997. 

\bibitem[Ch]{Ch} Y. {\sc Choi}. {\em Group representations with empty residual spectrum.}
To appear in  Int. Eq. Op. Th. 

\bibitem[CG]{CG} E. {\sc Cordero}, K. {\sc Gröchenig}, M. {\sc
Leinert}. \newblock {\em Localization of frames II.}  \newblock  Appl. Comput. Harmonic Anal.
17(2004), 29-47. 


\bibitem[CS]{CS} O. {\sc Christensen} and T. {\sc Strohmer}. {\em The finite section method and problems in frame theory}. \newblock  J. 
Approx. Th., 133(2005), 221-237. 


\bibitem[FGL1]{FGL} G. {\sc Fendler}, K. {\sc Gröchenig}, M. {\sc
Leinert}. \newblock {\em Symmetry of weighted $L^1$-algebras and
the GRS-condition.} Bull. London Math. Soc. 38 (4), 625-635, 2006.

\bibitem[FGL2]{FGL2} G. {\sc Fendler}, K. {\sc Gröchenig}, M. {\sc
Leinert}. \newblock {\em Convolution-Dominated Operators on Discrete Groups.} \newblock 
Integr. equ. oper. theory 61, 493-500, 2008. 



\bibitem[G1]{G1} K. {\sc Gröchenig}, M. {\sc
Leinert}. \newblock {\em Foundation of Time-Frequency Analysis}. Birkhäuser, Boston, 2001. 

\bibitem[G2]{G2} K. {\sc Gröchenig}. \newblock {\em Localized frames are finite unions of Riesz sequences}. Adv. Comput. Math., 
18(2003), 149-157. 


\bibitem[G3]{G3} K. {\sc Gröchenig}. \newblock {\em Localization of frames, Banach frames, and the invertibility of the frame op- 
erator}. J. Fourier Anal. Appl., 10(2004), 105-132. 

\bibitem[G4]{G4} K. {\sc Gröchenig}. \newblock {\em Time-frequency analysis of Sjöstrand's class}. Rev. Mat. Iberoam., To appear. 


\bibitem[GL]{GL} K. {\sc Gröchenig}, M. {\sc
Leinert}. \newblock {\em Symmetry of matrix algebras and symbolic 
calculus for infinite matrices}, Trans, Amer. Math. Soc., 358, 2695-2711, 2006. 


\bibitem[Gro1]{Gro} M. {\sc Gromov}.  \newblock {\em Groups of polynomial growth and expanding
maps.} \newblock Publ. Math. IHES, 53, 53-73, 1981.


\bibitem[Gui]{Gui} Y. {\sc Guivarc'h}.  \newblock {\em Croissance polynômiale et périodes des
fonctions harmoniques}. \newblock Bull. Sc. Math. France 101,
333-379, 1973.

\bibitem[J]{J} S. {\sc Jaffard}.  \newblock {\em Propriétés des matrices ``bien localisée" près de leur diagonale et quelques applications.} \newblock Ann. Inst. H. Poincaré Anal. Non lineaire, 7 (5), 461-473, 1990.



\bibitem[R]{R} J. {\sc Rosenblatt}, {\em Invariant measures and growth conditions}, Trans. Amer. Math. Soc. 193 (1974), 
p. 33-53. 

\bibitem[ShS]{Sun2} C. E. {\sc Shin} and Q. {\sc Sun}. \newblock  {\em Stability of localized operators}.  \newblock Journal of Functional Analysis,  256, 2417-2439, 2009.


\bibitem[S1]{S1} J. {\sc Stohmer}. \newblock {\em Pseudo-differential operators and Banach algebras in mobile communications.} \newblock Appl. Comput. Harmon. Anal., 20(2), 237-249, 2006.

\bibitem[S2]{S2} J. {\sc Stohmer}. \newblock {\em Rates of convergence for the approximation of shift-invariant systems in $\ell^2(\Z)$.} \newblock J. Fourier Anal. Appl., 5(2000), 519-616. 


\bibitem[S3]{S3} J. {\sc Stohmer}. \newblock {\em Four short stories about Toeplitz matrix calculations.} \newblock  Linear Algebra Appl., 
343/344(2002), 321-344. 

\bibitem[Sj]{Sj} J. {\sc Sjöstrand}. \newblock {\em Wiener type algebra of pseudodifferential operators, Centre de Mathematiques}, \newblock
Ecole Polytechnique, Palaiseau France, Seminaire 1994-1995, December 1994. 


\bibitem[Su]{Sun} Q. {\sc Sun}. \newblock {\em Wiener's lemma for infinite matrices.} \newblock Trans. Amer., Math. Soc., 359(7): 3099-3123 (electronic),  2007.


\bibitem[T]{T} R. {\sc Tessera}. \newblock {\em The inclusion of the Schur algebra in $B(\ell^2)$ is not inverse-closed.} \newblock Preprint, 2009.



\end{thebibliography}
\end{document}